\documentclass[11pt,leqno]{article}
\usepackage{amsmath}
\usepackage{amssymb}
\usepackage{amscd}

\renewcommand{\phi}{\varphi}

\newtheorem{theorem}{Theorem}[section]
\newtheorem{proposition}[theorem]{Proposition}
\newtheorem{lemma}[theorem]{Lemma}
\newtheorem{corollary}[theorem]{Corollary}

\newtheorem{definition}[theorem]{Definition}

\begin{document}
\thispagestyle{empty}

\begin{center}
{\Large\bf Extensions of Quasidiagonal $C^*$-algebras and K-theory}

\end{center}

\begin{center}{\bf Nathanial P. Brown \footnote{
      NSF Postdoctoral Fellow. }}\\
   UC-Berkeley\\
   Berkeley, California 94720 \\
   nbrown{\char'100}math.berkeley.edu

\vspace{2mm}

{\bf Marius Dadarlat\footnote{
      Partially supported by an NSF grant. }}

Purdue University\\
West Lafayette, Indiana 47907 \\
mdd{\char'100}math.purdue.edu

\end{center}

\begin{abstract}
Let $0 \to I \to E \to B \to 0$ be a short exact sequence of
C*-algebras where $E$ is separable, $I$ is quasidiagonal (QD) and $B$
is nuclear, QD and satisfies the UCT.  It is shown that if the
boundary map $\partial : K_1(B) \to K_0 (I)$ vanishes then $E$ must be
QD also.

A Hahn-Banach type property for $K_0$ of QD $C^*$-algebras is also
formulated.  It is shown that every nuclear QD $C^*$-algebra has this
$K_0$-Hahn-Banach property if and only if the boundary map $\partial :
K_1(B) \to K_0 (I)$ (from above) always completely determines when $E$
is QD in the nuclear case.
\end{abstract}

\parskip2mm

\section{Introduction}

Quasidiagonal (QD) $C^*$-algebras are those which enjoy a certain
finite dimensional approximation property. (See \cite{vo2}, \cite{Br3}
for surveys of the theory of QD $C^*$-algebras.) While these finite
dimensional approximations have certainly lead to a better
understanding of the structure of QD $C^*$-algebras, there are a
number of very basic open questions.  For example, assume that $0 \to
I \to E \stackrel{\pi}{\to} B \to 0$ is a split exact sequence (i.e.\
there exists a *-homomorphism $\phi:B \to E$ such that $\pi \circ \phi
= id_B$) where both $I$ and $B$ are QD.  It is not known whether $E$
must be QD (and, in fact, it is not even clear what to expect).

In this paper we study the extension problem for QD $C^*$-algebras and
it's relation to some natural questions concerning K-theory of QD
$C^*$-algebras.  Our techniques rely heavily on Kasparov's theory of
extensions and thus we will always need some nuclearity assumptions.

For example, adapting techniques found in \cite{Sp} we will show
(Theorem~\ref{3.4}) that if $0 \to I \to E \to B \to 0$ is short exact
where $E$ is separable, $I$ is QD, $B$ is nuclear, QD and satisfies
the Universal Coefficient Theorem (UCT) and the boundary map $\partial
: K_1(B) \to K_0 (I)$ vanishes then $E$ must be QD also.  It follows
that if $K_1 (B) = 0$ then $E$ is always QD, which generalizes work of
Eilers, Loring and Pedersen (\cite{ELP}).  As another application we
observe that in the case that $I$ is the compact operators our result
implies that $E$ is QD if and only if the (class of the) extension is
in the kernel of the natural map $Ext(B) \to Hom(K_1(B), {\mathbb
Z})$, where $Ext(B)$ denotes the classical BDF group (recall that we
are assuming $B$ is nuclear and hence $Ext(B)$ is a group).  Also, we
verify a conjecture of \cite{BK}, stating that an asymptotically split
extension of NF algebras is NF, under the additional assumption that
the quotient algebra satisfies the UCT of \cite{RS}.

We then study the general extension problem.  Now let $0 \to I \to E
\to B \to 0$ be exact where $E$ is separable and nuclear, $I$ is QD
and $B$ is QD and satisfies the UCT. Based on previous work of
Spielberg (\cite{Sp}) it is reasonable to expect that in this case $E$
will be QD if and only if $\partial (K_1(B)) \cap K_0^+ (I) = \{ 0
\}$, where $K_0^+ (I) = \{ 0 \}$ denotes the positive cone of $K_0
(I)$.  Though we are unable to resolve this question we do show that
it is equivalent to some other natural questions concerning the
K-theory of QD $C^*$-algebras and that in order to solve the general
extension problem it suffices to prove the special case that $B =
C({\mathbb T})$ (see Theorem~\ref{4.10}).

The first equivalent K-theory question is: If $A$ is nuclear,
separable and QD and $G \subset K_0 (A)$ is a subgroup such that $G
\cap K_0^+ (A) = 0$ then can one always find an embedding $\rho : A
\hookrightarrow C$ where $C$ is QD and $\rho_* (G) = 0$?  The
condition $G \cap K_0^+ (A) = 0$ is easily seen to be necessary and
hence the question is whether or not it is sufficient.  The second
K-theory question asks whether every nuclear QD $C^*$-algebra
satisfies what we call the {\em $K_0$-Hahn-Banach property} (see
Definition 4.7).  Roughly speaking this $K_0$-Hahn-Banach property
states that if $x \in K_0 (A)$ and $\pm x \notin K_0^+(A)$ then one
can always find finite dimensional approximate morphisms (i.e.\
"functionals") which separate $x$ from $K_0^+(A)$.  (Due to possible
perforation in $K_0(A)$ this statement is not quite correct, but it
conveys the main idea.)  Determining whether every nuclear QD algebra
satisfies the $K_0$-Hahn-Banach property is of independent interest as
our inability to understand how well finite dimensional approximate
morphisms read K-theory has been a major obstacle in the
classification program.

In section 2 we review the necessary theory of extensions and prove a
few simple results needed later.  In section 3 we handle the case when
$\partial : K_1(B) \to K_0 (I)$ vanishes.  In section 4 we turn to the
general extension problem and show equivalence with the K-theory
questions described above.

The present work is related to work of Salinas \cite{Sa},
\cite{Salinas} and Schochet \cite{Sch}. Those authors study the
quasidiagonality of extensions $0 \to I \to E \to B \to 0$ (i.e. the
question of whether or not $I$ contains an approximate unit of
projections which is quasicentral in $E$) whereas we study the QD of
the C*-algebra $E$. The two questions are different even if $I$ is the
compact operators. Indeed, while the quasidiagonality of $0 \to
\mathcal{K} \to E \to B \to 0$ does imply the QD of $E$, the converse
implication is false (see Section~\ref{3}).

\section{Preliminaries and Trivial Extensions.}

Most of this section is devoted to reviewing definitions, introducing
notation and recalling some standard facts about extensions of
$C^*$-algebras.  However, at the end we prove a few simple facts which
will be needed later.  The main result states that quasidiagonality is
preserved in split extensions provided that either the ideal or the
quotient is a nuclear $C^*$-algebra (see Proposition~\ref{2.5}).

For a comprehensive introduction to the aspects of extension theory
which we will need we recommend looking at \cite[Section 15]{Bl}. For
any $C^*$-algebra $I$ we will let $M(I)$ be it's multiplier algebra
and $Q(I) = M(I)/I$ be it's corona algebra.  Let $\pi : M(I) \to Q(I)$
be the quotient map.

If $E$ is any $C^*$-algebra containing $I$ as an ideal and $B = E/I$
then there exists a unique *-homomorphism $\rho : E \to M(I)$ such
that $\rho(I) = I$ and hence an induced *-homomorphism $\gamma : B \to
Q(I)$.  The map $\gamma$ is injective if and only if $\rho$ is in
injective if and only if $I$ sits as an essential ideal in $E$.
Conversely, given a $C^*$-algebra $B$ and a *-homomorphism $\gamma : B
\to Q(I)$ we can construct the pullback which, by definition, is the
$C^*$-algebra $$E(\gamma) = \{ x \oplus b \in M(I) \oplus B : \pi(x) =
\gamma(b) \}.$$ This gives a short exact sequence $0 \to I \to
E(\gamma) \to B \to 0$.  Moreover, if $B = E/I$ with induced map
$\gamma : B \to Q(I)$ then there is an induced *-isomorphism $\Phi : E
\to E(\gamma)$ with commutativity in the diagram
$$\begin{CD}
0 @>>> I @>>> E @>>> B @>>> 0 \\
@. @| @V \Phi VV @| @. \\
0 @>>> I @>>> E(\gamma) @>>> B @>>> 0.
\end{CD}$$
Hence there is a one to one correspondence between extensions of $I$
by $B$ and *-homomorphisms $\gamma : B \to Q(I)$.  As is standard, we
will refer to a *-homomorphism $\gamma : B \to Q(I)$ as a {\em Busby
invariant} and freely use the above correspondence between Busby
invariants and extensions.

When $I$ is stable (i.e.\ $I \cong {\cal K}\otimes I$, where ${\cal
K}$ denotes the compact operators on a separable infinite dimensional
Hilbert space) there is a natural way of adding two extensions which
we now describe.  Any isomorphism $M_2 ({\mathbb C}) \otimes {\cal K}
\cong {\cal K}$ induces an isomorphism $M_2 ({\mathbb C}) \otimes
{\cal K} \otimes I \cong {\cal K} \otimes I$ which then gives
isomorphisms $M_2 ({\mathbb C}) \otimes M({\cal K}\otimes I) \cong
M({\cal K}\otimes I)$ and $M_2 ({\mathbb C}) \otimes Q({\cal K}\otimes
I) \cong Q({\cal K}\otimes I)$.  Thus if we are given two Busby
invariants $\gamma_1, \ \gamma_2 : B \to Q({\cal K}\otimes I)$ we can
define a new Busby invariant $\gamma_1 \oplus \gamma_2$ by $$(\gamma_1
\oplus \gamma_2)(b) = \left(\begin{array}{cc} \gamma_1(b) & 0 \\ 0 &
\gamma_2(b)
\end{array}\right)
\in M_2 ({\mathbb C}) \otimes Q({\cal K}\otimes I) \cong Q({\cal K}\otimes I).$$

Of course the Busby invariant $\gamma_1 \oplus \gamma_2$ constructed
in this way will depend on the particular isomorphism $M_2 ({\mathbb
C}) \otimes {\cal K} \cong {\cal K}$.  To remedy this we say that two
Busby invariants $\gamma_1, \ \gamma_2$ are {\em strongly equivalent}
if there exists a unitary $u \in M(I)$ such that ${\rm Ad}\pi (u)\big(
\gamma_1 (b) \big) = \pi(u) \gamma_1 (b) \pi(u^*) = \gamma_2 (b)$, for
all $b \in B$, where $\pi : M(I) \to Q(I)$ is the quotient map.  Note
that if $\gamma_1$ and $\gamma_2$ are strongly equivalent then
$E(\gamma_1)$ and $E(\gamma_2)$ are isomorphic $C^*$-algebras. Indeed,
the map $E(\gamma_1) \to E(\gamma_2)$, $x \oplus b \mapsto uxu^*
\oplus b$ is easily seen to be an isomorphism.  Since any isomorphism
$M_2 ({\mathbb C}) \otimes {\cal K} \cong {\cal K}$ is implemented by
a unitary we see that $\gamma_1 \oplus \gamma_2$ is unique up to
strong equivalence.  In particular, the isomorphism class of the
$C^*$-algebra $E(\gamma_1 \oplus \gamma_2)$ does not depend on the
choice of isomorphism $M_2 ({\mathbb C}) \otimes {\cal K} \cong {\cal
K}$.

A Busby invariant $\tau$ is called {\em trivial} if it lifts to a
*-homomorphism $\phi : B \to M(I)$ (i.e.\ $\pi\circ \phi = \gamma$). A
Busby invariant $\gamma : B \to Q({\cal K}\otimes I)$ is called {\em
absorbing} if $\gamma \oplus \tau$ is strongly equivalent to $\gamma$
for every trivial $\tau$.  Note that if $\gamma$ is absorbing then so
is $\tilde{\gamma} \oplus \gamma$ for any $\tilde{\gamma}$. In
particular if $\gamma$ is absorbing then $\gamma$ is injective. Note
also that if $\tau_1$ and $\tau_2$ are both trivial and absorbing then
$\tau_1$, $\tau_1 \oplus \tau_2$ and $\tau_2$ are all strongly
equivalent.  Thus we get the following fact.

\begin{lemma}\label{2.1}
If $\tau_1, \ \tau_2 : B \to Q({\cal K}\otimes I)$ are both
trivial and absorbing then $E(\tau_1) \cong E(\tau_2)$.
\end{lemma}

Another simple fact we will need is the following.

\begin{lemma}\label{2.2}
If $\gamma, \ \tau : B \to Q({\cal K}\otimes I)$ are Busby invariants
with $\tau$ trivial then there is a natural embedding $E(\gamma)
\hookrightarrow E(\gamma \oplus \tau)$.
\end{lemma}

{\noindent\it Proof.}  Let $\phi : B \to M(I)$ be a lifting of
$\tau$. Define a map $E(\gamma) \to E(\gamma \oplus \tau)$ by $$x
\oplus b \mapsto
\left(\begin{array}{cc}
x  & 0 \\
0 & \phi(b)
\end{array}\right)
\oplus b.$$ Evidently this map is an injective $*$-homomorphism.  $\square$

The following generalization of Voiculescu's Theorem, which is due to
Kasparov, will be crucial in what follows.

\begin{theorem}\label{2.3}(\cite[Thm.~15.12.4]{Bl}) Assume that $B$ is
separable, $I$ is $\sigma$-unital and either $B$ or $I$ is nuclear.
Let $\rho : B \to B(H)$ be a faithful representation such that $H$ is
separable, $\rho(B) \cap {\cal K}(H) = \{ 0 \}$ and the orthogonal
complement of the nondegeneracy subspace of $\rho(B)$ (i.e.\ $H
\ominus \overline{\rho(B)H}$) is infinite dimensional.  Regarding
$B(H) \cong B(H)\otimes 1 \subset M({\cal K} \otimes I)$ as scalar
operators we get a short exact sequence $$0 \to {\cal K} \otimes I \to
\rho(B) \otimes 1 + {\cal K} \otimes I \to B \to 0.$$ If $\tau$
denotes the induced Busby invariant then $\tau$ is both trivial and
absorbing.
\end{theorem}

We define an equivalence relation on the set of Bubsy invariants $B
\to Q({\cal K} \otimes I)$ by saying $\gamma$ is related to
$\tilde{\gamma}$ if there exist trivial Busby invariants $\tau,
\tilde{\tau}$ such that $\gamma \oplus \tau$ is strongly equivalent to
$\tilde{\gamma} \oplus \tilde{\tau}$.  Taking the quotient by this
relation yields the semigroup $Ext(B, {\cal K} \otimes I)$.  The image
of a map $\gamma : B \to Q({\cal K} \otimes I)$ in $Ext(B, {\cal K}
\otimes I)$ is denoted $[\gamma]$.  Note that all trivial Busby
invariants give rise to the same class denoted by $0 \in Ext(B, {\cal
K} \otimes I)$ and this class is a neutral element (i.e.\ identity)
for the semigroup.  Note also that if $[\gamma] = 0 \in Ext(B, {\cal
K} \otimes I)$ then it does not follow that $\gamma$ is trivial.
However, it does follow that if $\tau$ is a trivial absorbing Busby
invariant then so is $\gamma \oplus \tau$.

We are almost ready to prove the main result of this section.  We just
need one more definition.

\begin{definition}\label{2.4}
If $0 \to I \to E \to B \to 0$ is an exact sequence with Busby
invariant $\gamma$ then we let $\gamma^s : {\cal K} \otimes B \to
Q({\cal K} \otimes I)$ denote the stabilization of $\gamma$.  That is,
$\gamma^s$ is the Busby invariant of the exact sequence $0 \to {\cal
K} \otimes I \to {\cal K} \otimes E \to {\cal K} \otimes B \to 0$.
\end{definition}

Note that there is always an embedding $E \cong E(\gamma)
\hookrightarrow E(\gamma^s)$.

\begin{proposition}\label{2.5}
Let $0 \to I \to E \to B \to 0$ be exact with Busby invariant
$\gamma$.  If both $I$ and $B$ are QD, $B$ is separable, $I$ is
$\sigma$-unital, either $I$ or $B$ is nuclear and $[\gamma^s] = 0 \in
Ext({\cal K} \otimes B, {\cal K} \otimes I)$ then $E$ is also QD.
\end{proposition}

{\noindent\it Proof.}  Since quasidiagonality passes to subalgebras,
it suffices to show that if $\tau : {\cal K} \otimes B \to Q({\cal K}
\otimes I)$ is a trivial absorbing Busby invariant (which exists by
Theorem~\ref{2.3}) then $E(\tau)$ is QD.  Indeed, by Lemmas~\ref{2.1},
\ref{2.2} and the definition of $Ext({\cal K} \otimes B, {\cal K}
\otimes I)$ we have the inclusions $$E \hookrightarrow E(\gamma^s)
\hookrightarrow E(\gamma^s \oplus \tau) \cong E(\tau).$$

To prove that $E(\tau)$ is QD we may assume (again by Lemma~\ref{2.1})
that $\tau$ arises from the particular extension described in
Theorem~\ref{2.3}. However for that extension it is easy to see that
$E(\tau) \hookrightarrow (\rho(B) + {\cal K}) \otimes \tilde{I}$,
where $\tilde{I}$ is the unitization of $I$.  But since $\rho(B) \cap
{\cal K} = \{ 0 \}$ it follows that $\rho(B) + {\cal K}$ is QD
(\cite[Thm.~3.11]{Br3}). Hence $(\rho(B) + {\cal K}) \otimes
\tilde{I}$ is also QD as a minimal tensor product QD C*-algebras
(\cite[Prop.~7.5]{Br3} ). $\square$

Note that the above proposition covers the case of split extensions
(i.e.\ when $\gamma$ is trivial).

\section{When $\partial : K_1 (B) \to K_0 (I)$ is zero.}\label{3}

The main result of this section (Theorem~\ref{3.4}) states that if the
boundary map $\partial : K_1 (B) \to K_0 (I)$ coming from an exact
sequence $0 \to I \to E \to B \to 0$ is zero then $E$ will be QD
whenever $I$ is QD and $B$ is nuclear, QD and satisfies the Universal
Coefficient Theorem (UCT) of Rosenberg and Schochet (\cite{RS}). The
main ideas in the proof are inspired by work of Spielberg
(\cite{Sp}). We also discuss a few consequences of our result,
including generalization of work of Eilers-Loring-Pedersen
(\cite{ELP}) and a partial solution to a conjecture of Blackadar and
Kirchberg \cite{BK}.

\begin{definition}\label{3.1}
{\em An embedding $I \hookrightarrow J$ is called {\em approximately
unital} if it takes an approximate unit of $I$ to an approximate unit
of $J$.  }
\end{definition}

In this case there is a natural inclusion $M(I) \hookrightarrow M(J)$
which induces an inclusion $Q(I) \hookrightarrow Q(J)$
\cite[3.12.12]{Pe1}.  Hence for any Busby invariant $\gamma : B \to
Q(I)$ there is an induced Busby invariant $\eta : B \to Q(J)$ with
commutativity in the diagram
$$\begin{CD}
0 @>>> I @>>> E(\gamma) @>>> B @>>> 0 \\
@. @VVV @VVV   @| @. \\
0 @>>> J @>>> E(\eta)  @>>> B @>>> 0.
\end{CD}$$
Moreover, the two vertical maps on the left are injective.

There are two ways of producing approximately unital embeddings which
we will need.  The first is $I \hookrightarrow I\otimes A$, for some
{\em unital} $C^*$-algebra $A$.  If $\{ e_{\lambda} \}$ is an
approximate unit of $I$ then, of course, $e_{\lambda} \otimes 1_A$
will be an approximate unit of $I \otimes A$. The other is to start
with an arbitrary embedding $I \hookrightarrow J^{\prime}$ and define
$J$ to be the hereditary subalgebra in $J^{\prime}$ generated by $I$.
That is, define $J$ to be the closure of $\cup_{\lambda} e_{\lambda}
J^{\prime} e_{\lambda}$.  One easily checks that $J$ is then a
hereditary subalgebra of $J^{\prime}$ and the embedding $I
\hookrightarrow J$ is approximately unital.

In the theory of separable QD $C^*$-algebras there are some {\em
nonseparable} algebras which play a key role.  The first is the direct
product $\Pi_i M_{n_i} ({\mathbb C})$ for some sequence of integers
$\{ n_i \}$.  This algebra is the multiplier algebra of the direct sum
$\oplus_i M_{n_i} ({\mathbb C})$.  If $H$ is any separable Hilbert
space then we can always find a decomposition $H = \oplus_i {\mathbb
C}^{n_i}$ and then we have natural inclusions $\oplus_i M_{n_i}
({\mathbb C}) \hookrightarrow {\cal K}(H)$, $\Pi_i M_{n_i} ({\mathbb
C}) \hookrightarrow B(H)$ and $Q(\oplus_i M_{n_i} ({\mathbb C}))
\hookrightarrow Q({\cal K}(H))$.  Another algebra which we will need
is $\Pi_i M_{n_i} ({\mathbb C}) + {\cal K}(H)$.

\begin{lemma}\label{3.2}
Let $J \subset \Pi_i M_{n_i} ({\mathbb C}) + {\cal K}(H)$ be a
hereditary subalgebra containing ${\cal K}(H)$.  Then $K_1 (J) = 0$.
\end{lemma}

{\noindent\it Proof.}  Letting $\pi : B(H) \to Q(H)$ be the quotient
map we have that $\pi(J)$ is a hereditary subalgebra of $Q(\oplus_i
M_{n_i} ({\mathbb C}))$ (use the fact that if $0 \leq a \in J, b \in
Q(\oplus_i M_{n_i} ({\mathbb C}))$ and $0 \leq b \leq \pi(a)$ then
there exists $0 \leq c \in \Pi_i M_{n_i} ({\mathbb C}) + {\cal K}(H)$
such that $c \leq a$ and $\pi(c) = b$; \cite[Cor.\ IX.4.5]{Da}.  Also,
the exact sequence $ 0 \to {\cal K}(H) \to J \to \pi(J) \to 0$ is a
{\em quasidiagonal extension} (i.e.\ ${\cal K}(H)$ contains an
approximate unit {\em of projections which is quasicentral in} $J$).
Hence \cite[Thm.~8]{BD}, states that we have a short exact sequence
$$0 \to K_1 ( {\cal K}(H)) \to K_1 ( J) \to K_1 ( \pi(J)) \to 0.$$
Thus it suffices to show that $K_1 (X) = 0$ for any hereditary
subalgebra $X$ of $Q(\oplus_i M_{n_i} ({\mathbb C}))$.

But if $X \subset Q(\oplus_i M_{n_i} ({\mathbb C}))$ is a hereditary
subalgebra then we can find a quasidiagonal extension $$0 \to \oplus_i
M_{n_i} ({\mathbb C}) \to Y \to X \to 0,$$ where $Y \subset \Pi_i
M_{n_i} ({\mathbb C})$ is a hereditary subalgebra.  Applying
\cite[Thm.\ 8]{BD} again it suffices to show that every hereditary
subalgebra of $\Pi_i M_{n_i} ({\mathbb C})$ has trivial $K_1$-group.

But, if $Y \subset \Pi_i M_{n_i} ({\mathbb C})$ is a hereditary
$\sigma$-unital subalgebra then $Y$ has an increasing approximate unit
consisting of projections, say $\{ e_{n} \}$ (\cite{BP}). Hence $$K_1
(Y) = \lim K_1 (e_{n} \Pi_i M_{n_i} ({\mathbb C}) e_{n}),$$ since $ Y
= \lim e_{n} \Pi_i M_{n_i} ({\mathbb C}) e_{n}$ (by heredity).  But
for each $n$ it is clear that $e_{n} \Pi_i M_{n_i} ({\mathbb C})
e_{n}$ is isomorphic to $\Pi_i M_{k_i} ({\mathbb C})$ for some
integers $\{ k_i \}$ and hence $K_1 (e_{n} \Pi_i M_{n_i} ({\mathbb C})
e_{n}) = 0$. $\square$

\begin{proposition}\label{3.3}
Let $I$ be a separable QD $C^*$-algebra.  Then there exists an
approximately unital embedding $I \hookrightarrow J$, where $J$ is a
 $\sigma$-unital QD $C^*$-algebra with $K_1 (J) =
0$.
\end{proposition}

{\noindent\it Proof.}  Let $\rho : I \to B(H)$ be a nondegenerate
faithful representation such that $\rho(I) \cap {\cal K}(H) = \{ 0
\}$.  By \cite[Prop.~5.2]{Br3}, there exists a decomposition $H =
\oplus_i {\mathbb C}^{n_i}$ such that $\rho(I) \subset \Pi_i M_{n_i}
({\mathbb C}) + {\cal K}(H)$.  Let $J$ be the hereditary subalgebra of
$\Pi_i M_{n_i} ({\mathbb C}) + {\cal K}(H)$ generated by
$\rho(I)$. The conclusion follows from the previous lemma. $\square$

For the remainder of this section we will let ${\cal U} = \otimes_n
M_n ({\mathbb C})$ be the Universal UHF algebra (i.e.\ the UHF algebra
with $K_0 ({\cal U}) = {\mathbb Q}$).  For any Busby invariant $\gamma
: B \to Q(J)$ we let $\gamma^{{\mathbb Q}}$ denote the Busby invariant
coming from the short exact sequence $$ 0 \to J\otimes {\cal U} \to
E(\gamma) \otimes {\cal U} \to B \otimes {\cal U} \to 0.$$

\begin{theorem}\label{3.4}
Let $0 \to I \to E \to B \to 0$ be a short exact sequence where $E$ is
separable, $I$ is QD and $B$ is nuclear, QD and satisfies the UCT.  If
the induced map $\partial : K_1 (B) \to K_0 (I)$ is zero then $E$ is
QD.
\end{theorem}

{\noindent\it Proof.}  Let $\gamma$ be the Busby invariant of the
exact sequence in question.  By the previous proposition we can find
an approximately unital embedding $I \hookrightarrow J$, where $J$ is
QD with $K_1 (J) = 0$.  By the remarks following Definition 3.1 we
have an inclusion $E \hookrightarrow E(\eta)$ where $\eta : B \to
Q(J)$ is the induced Busby invariant.  By naturality we then have that
both index maps $ \partial : K_1 (B) \to K_0 (J)$ and $\partial : K_0
(B) \to K_1 (J)$ are zero.  Hence the index maps arising from the
stabilization $\eta^s : B\otimes {\cal K} \to Q(J \otimes {\cal K})$
are also zero.

Now, if it happens that $K_0 (J)$ is a divisible group then the
Universal Coefficient Theorem would imply that $[\eta^s] = 0 \in
Ext(B\otimes {\cal K}, J \otimes {\cal K})$ and so by
Proposition~\ref{2.5} we would be done.  Of course this will not be
true in general and so may have to replace $\eta^s$ with
$(\eta^s)^{{\mathbb Q}}$.  But applying naturality one more time, both
boundary maps on K-theory arising from $(\eta^s)^{{\mathbb Q}}$ will
also vanish.  Hence the theorem follows from the inclusions $E
\hookrightarrow E(\eta) \hookrightarrow E(\eta^s) \hookrightarrow
E((\eta^s)^{{\mathbb Q}})$ together with Proposition~\ref{2.5} applied
to $(\eta^s)^{{\mathbb Q}}$.  $\square$

In the case that the ideal is nuclear and the quotient is an AF
algebra, the next result was obtained by Eilers, Loring and Pedersen
(\cite[ Cor.~4.6]{ELP}).

\begin{corollary}\label{3.5}
Assume that $B$ is a separable nuclear QD $C^*$-algebra satisfying the
UCT and with $K_1 (B) = 0$.  For any separable QD $C^*$-algebra $I$
and Busby invariant $\gamma : B \to Q(I)$ we have that $E(\gamma)$ is
QD.
\end{corollary}

This corollary actually extends to the case where $K_1 (B)$ is a
torsion group since we can tensor any short exact sequence with ${\cal
U}$ and $K_1 (B \otimes {\cal U}) = 0$ in this case.  For example,
this would cover the case that $B = C_0(\mathbb{R}) \otimes {\cal
O}_n, \ (2 \leq n \leq \infty),$ where ${\cal O}_n$ denotes the Cuntz
algebra on $n$ generators.  Similarly, it is clear that
Theorem~\ref{3.4} is valid under the weaker hypothesis that $\partial
(K_1 (B))$ is contained in the torsion subgroup of $K_0 (I)$.

\begin{definition}\label{3.6}
{\em For any two QD $C^*$-algebras $I, \ B$ let $Ext_{QD} (B,{\cal K}
\otimes I) \subset Ext(B,{\cal K} \otimes I)$ denote the set of
classes of Busby invariants $\gamma$ such that $E(\gamma)$ is QD.  }
\end{definition}

It is easy to check that if $[\gamma] = [\tilde{\gamma}] \in
Ext(B,{\cal K} \otimes I)$ then $E(\gamma)$ is QD if and only if
$E(\tilde{\gamma})$ is QD and hence $Ext_{QD} (B,{\cal K} \otimes I)$
is well defined.  It is also easy to see that $Ext_{QD} (B,{\cal K}
\otimes I)$ is a sub-semigroup of $Ext(B,{\cal K} \otimes I)$.
Finally, we remark that in the case $I = {\mathbb C}$ we do {\em not}
get the semigroup $Ext_{qd} (B,{\cal K})$ defined by Salinas; it
follows from Corollary~\ref{3.7} below, however, that we do get what
he called $Ext_{bqt} (B,{\cal K})$ in this case (see \cite[Definitions
2.7, 2.12 and Thm.~2.14]{Sa}). One has $Ext_{qd} (B,{\cal K}) \subset
Ext_{QD} (B,{\cal K} )$. The elements of $Ext_{QD} (B,{\cal K})$
corresponds to C*-algebras $E(\gamma)$ that are QD whereas $[\gamma]
\in Ext_{qd} (B,{\cal K})$ if the only if the extension $0 \to
\mathcal{K} \to E(\gamma) \to B \to 0$ is QD i.e. the concrete set
$E(\gamma) \subset M(\cal K)$ is QD.

Recall that there is a natural group homomorphism $\Phi : Ext(B, {\cal
K} \otimes I) \to Hom(K_1(B), K_0 (I))$ taking a Busby invariant to
the corresponding boundary map on K-theory.  From Theorem~\ref{3.4} it
follows that we always have an inclusion $Ker(\Phi) \subset Ext_{QD}
(B,{\cal K} \otimes I)$, when $B$ is nuclear, QD and satisfies the
UCT.  In general this inclusion will be proper, but we now describe a
class of algebras for which we have equality.

There is a natural semigroup $K_0^+ (I) \subset K_0 (I)$, called the
{\em positive cone}, given by $$K_0^+ (I) = \bigcup_{n \in {\mathbb
N}} \{ x \in K_0 (I) : x = [p], \ {\rm for \ some \ projection} \ p
\in M_n (I) \}.$$ When $I$ is unital this semigroup generates $K_0(I)$
but can also be trivial in general (e.g.\ if $I$ is stably
projectionless).  The natural isomorphism $K_0 (I) \cong K_0 ({\cal K}
\otimes I)$ induced by an embedding $I = e_{11} \otimes I \subset
{\cal K} \otimes I$, where $e_{11}$ is a minimal projection in ${\cal
K}$, preserves the positive cones.  We say that $K_0 (I)$ is {\em
totally ordered} if for every $x \in K_0 (I)$ either $x$ or $-x$ is an
element of $K_0^+ (I)$.

\begin{corollary}\label{3.7}
Assume $I$ is separable, QD and $K_0(I)$ is totally ordered.  For any
separable, nuclear, QD algebra $B$ which satisfies the UCT we have
that $Ext_{QD}(B, {\cal K} \otimes I) = Ker(\Phi)$.
\end{corollary}

{\noindent\it Proof.}  We only have to show $Ext_{QD} (B,{\cal K}
\otimes I) \subset Ker(\Phi)$. So let $[\gamma] \in Ext (B,{\cal K}
\otimes I)$.  If $E(\gamma)$ is a stably finite $C^*$-algebra then a
result of Spielberg (see Proposition~\ref{4.1} of the next section),
together with the assumption that $K_0 (I)$ is totally ordered,
implies that $[\gamma] \in Ker(\Phi)$.  But since QD implies stably
finite (\cite[Prop.~3.19]{Br3}) we have that if $[\gamma] \in Ext_{QD}
(B,{\cal K} \otimes I)$ then $[\gamma] \in Ker(\Phi)$. $\square$

The classic example for which $K_0(I)$ is totally ordered is the case
when $I = {\cal K}$.  In this setting the corollary above is very
similar to a result of Salinas' which describes the closure of $0 \in
Ext(B, {\cal K})$ in terms of quasidiagonality
(\cite[Thm.~2.9]{Sa}). See also \cite[Thm.~2.14]{Sa} for another
characterization of $Ext_{QD} (B, {\cal K})$ in terms of
bi-quasitriangular operators. For a K-theoretical characterization of
$Ext_{qd} (B, {\cal K})$ see \cite[Theorem 8.3]{Sch}.

The class of NF algebras introduced in \cite{BK} coincides with the
class of separable QD nuclear C*-algebras. It was conjectured in
\cite[Conj. 7.1.6]{BK} that an asymptotically split extension of NF
algebras is NF. We can verify the conjecture under an additional
asumption.

\begin{corollary}\label{3.8}
Let $0 \to I \to E \to B \to 0$ be an asymptotically split extension
with $I$ and $B$ NF algebras. If $B$ satisfies the UCT, then $E$ is
NF.
\end{corollary}

{\noindent\it Proof.} Both index maps are vanishing since the extension is
asymptotically split. The conclusion follows from Theorem~\ref{3.4}. $\square$

\section{Extensions and K-theory}

In this section we show that the general extension problem for nuclear
QD $C^*$-algebras is equivalent to some natural K-theoretic questions.

We begin by recalling a result of Spielberg which solves the extension
problem for stably finite $C^*$-algebras and shows that it is
completely governed by K-theory.

\begin{proposition}\label{4.1} \cite[Lemma 1.5]{Sp}
Let $0 \to I \to E \to B \to 0$ be short exact where both $I$ and $B$
are stably finite.  Then $E$ is stably finite if and only if
$\partial(K_1 (B)) \cap K_0^+ (I) = \{ 0 \},$ where $\partial : K_1
(B) \to K_0 (I)$ is the boundary map of the sequence.
\end{proposition}

In \cite[Question 7.3.1]{BK}, it is asked whether every nuclear stably
finite $C^*$-algebra is QD.  Support for an affirmative answer to this
question is provided by a number of nontrivial examples (\cite{Pi},
\cite{Sp}, \cite{Br1}, \cite{Br2}).  In fact, Corollary~\ref{3.7}
above also provides examples since the proof shows the equivalence of
quasidiagonality and stable finiteness (in fact we did not even assume
nuclearity of $E$ in that corollary).  Hence it is natural to wonder
if Spielberg's criterion completely determines quasidiagonality in
extensions as well. The following result gives some more evidence for
an affirmative answer. If $I$ is a C*-algebra, let
$SI=C_0(\mathbb{R})\otimes I$ denote the suspension of $I$. Note that
$K_0(SI)^{+}=\{0\}$ since $SI \otimes \mathcal{K}$ contains no nonzero
projections.

\begin{proposition} Let $0 \to SI \to E \to B \to 0$
be exact, where $I$ is $\sigma$-unital and $B$ is separable, QD,
nuclear. Then $E$ is QD.
\end{proposition}

{\noindent\it Proof.} The suspension $SI$ of $I$ is QD by
\cite{vo1}. We may assume that $I$ is stable. Let $\alpha:SI
\hookrightarrow SI$ be a null-homotopic approximately unital embedding
and let $\widehat{\alpha}:Q(SI)\hookrightarrow Q(SI)$ be the
corresponding $*$-monomorphism. Then for any Busby invariant $\gamma:B
\to M(SI)$, $[\widehat{\alpha}\circ \gamma]=0 \in Ext(B,SI)$ by the
homotopy invariance of $Ext(B,SI)$ in the second variable
\cite{Kas}. It follows that $E(\gamma) \hookrightarrow
E(\widehat{\alpha}\circ \gamma)$ is QD by Proposition~\ref{2.5}.
$\square$

\begin{definition}\label{4.2}
{\em Say that a QD $C^*$-algebra $A$ has the {\em QD extension
property} if for every separable, nuclear, QD algebra $B$ which
satisfies the UCT and Busby invariant $\gamma : B \to Q({\cal K}
\otimes A)$ we have that $E(\gamma)$ is QD if and only if $E(\gamma)$
is stably finite (which is if and only if $\partial (K_1(B)) \cap
K_0^+ ({\cal K} \otimes A) = \{ 0 \}$, by Proposition~\ref{4.1}).  }
\end{definition}

The QD extension property is closely related to a certain embedding
property for the K-theory of $A$ which we now describe.  The interest
in controlling the K-theory of embeddings of $C^*$-algebras goes back
to the seminal work of Pimsner and Voiculescu on AF embeddings of
irrational rotation algebras (\cite{PV}).  Since then other authors
have studied the K-theory of (AF) embeddings (\cite{Lo}, \cite{EL},
\cite{DL}, \cite{Br1}, \cite{Br1}).

\begin{definition}\label{4.3}
{\em Say that a QD $C^*$-algebra $A$ has the {\em $K_0$-embedding property}
if for every subgroup $G \subset K_0 (A)$ such that $G \cap K_0^+ (A)
= \{ 0 \}$ there exists an embedding $\rho : A \hookrightarrow C$,
where $C$ is also QD, such that $\rho_* (G) = 0$.}
\end{definition}

It is not hard to see that if $C$ is a stably finite $C^*$-algebra and
$p \in C$ is a nonzero projection then $[p]$ must be a nonzero element
of $K_0 (C)$.  From this remark it follows that the condition $G \cap
K_0^+ (A) = \{ 0 \}$ is necessary.  Hence the $K_0$-embedding property
states that this condition is also sufficient.

A number of QD $C^*$-algebras have the $K_0$-embedding property.  For
example, commutative C*-algebras, AF algebras (\cite[Lem.~1.14]{Sp}),
crossed products of AF algebras by ${\mathbb Z}$
(\cite[Thm.~5.5]{Br1}) and simple nuclear unital C*-algebras with
unique trace.

Our next goal is to connect the QD extension and $K_0$-embedding
properties.  But we first need a simple lemma.

\begin{lemma}\label{4.4}
Let $C$ be a hereditary subalgebra of a unital C*-algebra $D$. If $C$
has an approximate unit consisting of projections and $K_0 (D)$ has
cancellation then the inclusion $C \hookrightarrow D$ induces an
injective map $K_0 (C) \hookrightarrow K_0 (D)$.
\end{lemma}

{\noindent\it Proof.}  By \emph{cancellation} we mean that if $p,q \in
M_n (D)$ are projections with $[p] = [q]$ in $K_0(D)$ then there
exists a partial isometry $v\in M_n(D)$ such that $v v^* = p$ and $v^*
v = q$.

Let $x = [p] - [q] \in K_0 (C)$ be an element such that $x = 0 \in K_0
(D)$.  Since $C$ has an approximate unit of projections, say $\{
e_{\lambda} \}$, we may assume that $p$ and $q$ are projections in
$(e_{\lambda} \otimes 1) C \otimes M_n ({\mathbb C}) (e_{\lambda}
\otimes 1)$ for sufficiently large $n$ and $\lambda$.  Since $[p] =
[q]$ in $K_0 (D)$ and this group has cancellation we can find a
partial isometry $v \in M_n (D)$ such that $v v^* = p$ and $v^* v =
q$.

We claim that actually $v \in M_n (C)$ (which will evidently prove the
lemma).  To see this we first note that $v = v v^* (v) v^* v = pvq$
and hence $$v = pvq = (e_{\lambda} \otimes 1)pvq(e_{\lambda} \otimes
1) = (e_{\lambda} \otimes 1) v (e_{\lambda} \otimes 1).$$ Hence $v \in
(e_{\lambda} \otimes 1) D \otimes M_n ({\mathbb C}) (e_{\lambda}
\otimes 1)$.  But since $C$ is hereditary in $D$, $C \otimes M_n
({\mathbb C}) $ is hereditary in $D \otimes M_n ({\mathbb C}) $ and
thus $$v \in (e_{\lambda} \otimes 1) D \otimes M_n ({\mathbb C})
(e_{\lambda} \otimes 1) \subset C \otimes M_n ({\mathbb C}).  \ \ \ \
\ \square$$

\begin{proposition}\label{4.5}
Let $A$ be a separable QD $C^*$-algebra.  Then $A$ satisfies the QD
extension property if and only if $A$ satisfies the $K_0$-embedding
property.
\end{proposition}

{\noindent\it Proof.}  We begin with the easy direction.  Assume that
$A$ has the QD extension property and let $G \subset K_0 (A)$ be a
subgroup such that $G \cap K_0^+ (A) = \{ 0 \}$.  Since abelian
$C^*$-algebras satisfy the UCT we can construct an extension $$0 \to
{\cal K}\otimes A \to E \to \oplus_{\mathbb N} C({\mathbb T}) \to 0,$$
such that $\partial(K_1 (\oplus_{\mathbb N} C({\mathbb T}))) =
\partial( \oplus_{\mathbb N} {\mathbb Z}) = G$.  But since $A$ has the
QD extension property $E$ must be a QD $C^*$-algebra.  Thus the
six-term K-theory exact sequence implies that $A$ has the
$K_0$-embedding property (i.e.\ the embedding into $E$ will work).

Conversely, assume that $A$ has the $K_0$-embedding property and let
$$0 \to {\cal K}\otimes A \to E \to B \to 0$$ be a short exact
sequence where $B$ is separable, nuclear, QD, satisfies the UCT and
$E$ is stably finite.

Let $G = \partial(K_1 (B)) \subset K({\cal K}\otimes A) \cong K_0
(A)$.  Since $E$ is stably finite, $G \cap K_0^+ (A) = \{ 0 \}$.  By
the $K_0$-embedding property we can find a QD $C^*$-algebra $C$ and an
embedding $\rho : A \hookrightarrow C$ such that $\rho_* (G) = 0$.
Since $A$ is separable we may assume that $C$ is also separable.
Indeed $K_0 (A)$ (and hence $G$) is countable.  Thus it only takes a
countable number of projections and partial isometries in matrices
over $C$ to kill off $\rho_* (G)$.  From this observation it is easy
to see that we may assume that $C$ is also separable.

Let $\pi : C \hookrightarrow \Pi_i M_{n_i} ({\mathbb C}) + {\cal
K}(H)$ be an embedding (the existence of which is ensured by the
separability of $C$) as in the proof of Proposition~\ref{3.3}.  Let $J
\subset \Pi_i M_{n_i} ({\mathbb C}) + {\cal K}(H)$ be the hereditary
subalgebra generated by $\pi \circ \rho (A)$.  Since $\Pi_i M_{n_i}
({\mathbb C}) + {\cal K}(H)$ has real rank zero and stable rank one it
follows from Lemma~\ref{4.4} that the inclusion $J \hookrightarrow
\Pi_i M_{n_i} ({\mathbb C}) + {\cal K}(H)$ induces an injective map
$K_0 (J) \hookrightarrow K_0 ( \Pi_i M_{n_i} ({\mathbb C}) + {\cal
K}(H))$.  Since $G$ is in the kernel of the K-theory map induced by
the embedding $\pi \circ \rho : A \to \Pi_i M_{n_i} ({\mathbb C}) +
{\cal K}(H)$ it follows that $G$ is also in the kernel of the K-theory
map induced by the embedding $\pi \circ \rho : A \to J$.  But the
embedding into $J$ is approximately unital by construction and so we
get a commutative diagram
$$\begin{CD}
0 @>>> {\cal K}\otimes A @>>> E @>>> B @>>> 0 \\
@. @VVV @VVV   @| @. \\
0 @>>> {\cal K}\otimes J @>>> E(\eta)  @>>> B @>>> 0,
\end{CD}$$
where $\eta$ is the induced Busby invariant and the two vertical
maps on the left are injective.

Now we are done since naturality of the boundary map implies that the
homomorphism $\partial : K_1(B) \to K_0 ({\cal K}\otimes J)$ is zero
and hence $E(\eta)$ is QD by Theorem~\ref{3.4}.  $\square$

We now wish to point out a connection between extensions of QD
$C^*$-algebras and another very natural K-theoretic question.  For
brevity, we say a linear map $\phi : A \to B$ is {\em ccp} if it is
contractive and completely positive (\cite{Pa}). We recall a theorem
of Voiculescu.

\begin{theorem}\label{4.6} \cite[Thm.~1]{vo1}
Let $A$ be a separable $C^*$-algebra.  Then $A$ is QD if and only if
there exists an asymptotically multiplicative, asymptotically
isometric sequence of ccp maps $\phi_n : A \to M_{k_n} ({\mathbb C})$
for some sequence of natural numbers $k_n$ (i.e.\ $\| \phi_n (ab) -
\phi_n (a) \phi_n (b) \| \to 0$ and $\| \phi_n (a) \| \to \| a \|$ for
all $a,b \in A$).
\end{theorem}

Given this abstract characterization of QD $C^*$-algebras it is
natural to ask how well these approximating maps  capture the
relevant K-theoretic data.

\begin{definition}\label{4.7}
{\em Say that a QD $C^*$-algebra $A$ has the {\em $K_0$-Hahn-Banach
property} if for each $x \in K_0 (A)$ such that ${\mathbb Z}x \cap
K_0^+ (A) = \{ 0 \}$, where ${\mathbb Z}x = \{ kx : k \in {\mathbb Z}
\}$, there exists a sequence of asymptotically multiplicative,
asymptotically isometric ccp maps $\phi_n : A \to M_{k_n} ({\mathbb
C})$ such that $(\phi_n)_* (x) = 0$ for all $n$ large enough.  }
\end{definition}

It is easy to see that if $y \in K_0 (A)$ and there exists a nonzero
 integer $k$ such that $ky \in K_0^+ (A)$ then {\em for every}
 asymptotically multiplicative, asymptotically isometric sequence of
 ccp maps $\phi_n : A \to M_{k_n} ({\mathbb C})$ we have $(\phi_n)_*
 (y) > 0$ (if $k > 0$) or $(\phi_n)_* (y) < 0$ (if $k < 0$), for all
 sufficiently large $n$.  Hence this $K_0$-Hahn-Banach property states
 that one can separate elements $x \in K_0 (A)$ such that ${\mathbb
 Z}x \cap K_0^+ (A) = \{ 0 \}$ from (finite subsets of) the positive
 cone using finite dimensional approximate morphisms.

Another way of thinking about this property is that $A$ has the
$K_0$-Hahn-Banach property if and only if finite dimensional
approximate morphisms determine the order on $K_0 (A)$ to a large
extent.  A more precise formulation is contained in the next
proposition (not needed for the rest of the paper).

\begin{proposition}\label{4.8}
The $K_0$-Hahn-Banach property is equivalent to the following
property: If $x \in K_0 (A)$ and for every sequence of asymptotically
multiplicative, asymptotically isometric ccp maps $\phi_n : A \to
M_{k_n} ({\mathbb C})$ we have that $(\phi_n)_* (x) > 0$ for all large
$n$ then there exists a positive integer $k$ such that $kx \in K_0^+
(A)$.
\end{proposition}

{\noindent\it Proof.}  We first show that the (contrapositive of the)
second property above follows from the $K_0$-Hahn-Banach property.  So
assume we are given an element $x \in K_0 (A)$ and assume that there
is {\em no} positive integer $k$ such that $kx \in K_0^+ (A)$.  We
must exhibit a sequence of asymptotically multiplicative,
asymptotically isometric ccp maps $\phi_n : A \to M_{k_n} ({\mathbb
C})$ such that $(\phi_n)_* (x) \leq 0$ for all sufficiently large
$n$. There are two cases.

If there exists a negative integer $k$ such that $kx \in K_0^+ (A)$
then for every sequence $\phi_n : A \to M_{k_n} ({\mathbb C})$ we have
$(\phi_n)_* (x) < 0$ for all sufficiently large $n$ (see the
discussion following definition 4.7).  The second case is if ${\mathbb
Z}x \cap K_0^+ (A) = \{ 0 \}$.  This case is obviously handled by the
$K_0$-Hahn-Banach property.

Now we show how the second property above implies the
$K_0$-Hahn-Banach property.  So let $x \in K_0 (A)$ be such that
${\mathbb Z}x \cap K_0^+ (A) = \{ 0 \}$.  Since no positive multiple
of $x$ is in $K_0^+ (A)$ the second property implies that we can find
some sequence $\phi_n : A \to M_{k_n} ({\mathbb C})$ such that
$(\phi_n)_* (x) \leq 0$ for all sufficiently large $n$.  Similarly,
since no positive multiple of $-x$ is in $K_0^+ (A)$ we can find a
sequence $\psi_n : A \to M_{j_n} ({\mathbb C})$ such that $(\psi_n)_*
(x) \geq 0$ for all sufficiently large $n$.  If either of $\{ \phi_n
\}$ or $\{ \psi_n \}$ contains a subsequence with equality at $0$ then
we are done so we assume that $(\phi_n)_* (x) = -s_n < 0$ and
$(\psi_n)_* (x) = t_n > 0$ for all (sufficiently large) $n$.  It is
now clear what to do: we simply add up appropriate numbers of copies
of $\phi_n$ and $\psi_n$ so that these positive and negative ranks
cancel.  More precisely we define maps $$\Phi_n =
(\bigoplus\limits_{1}^{t_n} \phi_n ) \oplus
(\bigoplus\limits_{1}^{s_n} \psi_n )$$ and regard these maps as taking
values in the $(t_n k_n + s_n j_n) \times (t_n k_n + s_n j_n)$
matrices.  $\square$

\begin{proposition}\label{4.9}
If a separable QD $C^*$-algebra $A$ has the QD extension property or,
equivalently, the $K_0$-embedding property then $A$ also has the
$K_0$-Hahn-Banach property.
\end{proposition}

{\noindent\it Proof.} Assume that $A$ has the $K_0$-embedding property
and we are given $x \in K_0 (A)$ such that ${\mathbb Z}x \cap K_0^+
(A) = \{ 0 \}$, where ${\mathbb Z}x = \{ kx : k \in {\mathbb Z} \}$.
By the $K_0$-embedding property we can find an embedding $\rho : A
\hookrightarrow C$, where $C$ is QD and $\rho_* (x) = 0$.  As in the
proof of Proposition~\ref{4.5} we may assume that $C$ is also
separable. But then take any asymptotically multiplicative,
asymptotically isometric sequence of contractive completely positive
maps $\phi_n : C \to M_{k_n} ({\mathbb C})$ and we get that $(\phi_n
\circ \rho)_* (x) = 0$ for all sufficiently large $n$.  $\square$

We do not know if the converse of the previous proposition holds.
However our final result will complete the circle for the class of
nuclear $C^*$-algebras.  Moreover, the next theorem also states that in
order to prove that every separable, nuclear, QD $C^*$-algebra has any
of the properties we have been studying, it actually suffices to
consider very special cases of either the QD extension problem or
$K_0$-embedding problem.

\begin{theorem}\label{4.10}
The following statements are equivalent.
\begin{enumerate}
\item Every separable, nuclear, QD $C^*$-algebra has
the QD extension property.
\item Every separable, nuclear, QD $C^*$-algebra has
the $K_0$-embedding property.
\item Every separable, nuclear, QD $C^*$-algebra has
the $K_0$-Hahn-Banach property.
\item If $A$ is any separable, nuclear, QD $C^*$-algebra and $x \in
K_0(A)$ is such that ${\mathbb Z}x \cap K_0^+ (A) = \{ 0 \}$ then
there exists an embedding $\rho : A \hookrightarrow C$, where $C$ is
QD (but not necessarily separable or nuclear), such that $\rho_* (x) =
0$.
\item If $A$ is any separable, nuclear, QD $C^*$-algebra and $x \in
K_0(A)$ is such that ${\mathbb Z}x \cap K_0^+ (A) = \{ 0 \}$ then
there exists a short exact sequence $0 \to {\cal K}\otimes A \to E \to
C({\mathbb T}) \to 0$ where $E$ is QD and $x \in \partial(K_1
(C({\mathbb T}))) = \partial( {\mathbb Z})$.
\end{enumerate}
\end{theorem}

{\noindent\it Proof.}  The proof of Proposition~\ref{4.5} carries over
verbatim to show the equivalence of 1 and 2.  That proof also shows
the equivalence of 4 and 5.  The previous proposition shows that 2
implies 3 and hence we are left to show that 3 implies 5 and 4 implies
2.

We begin with the easier implication 4 $\implies$ 2.  So, let $A$ be
any separable, nuclear, QD $C^*$-algebra and $G \subset K_0 (A)$ be a
subgroup such that $G \cap K_0^+ (A) = \{ 0 \}$.  As in the proof of
Proposition~\ref{4.5} we can construct a short exact sequence $$0 \to
{\cal K}\otimes A \to E \to \bigoplus\limits_{1}^{\infty} C({\mathbb
T}) \to 0,$$ such that $\partial(K_1 (\oplus_{\mathbb N} C({\mathbb
T}))) = \partial( \oplus_{\mathbb N} {\mathbb Z}) = G$.  We will prove
that $E$ is QD and, by exactness of $\oplus_{\mathbb N} {\mathbb
Z}\stackrel{\partial}{\to} K_0(A) \to K_1(E)$, this will show 2.

For each $n$ there is a short exact sequence $$0 \to {\cal K}\otimes A
\to E_n \to \bigoplus\limits_{1}^{n} C({\mathbb T}) \to 0,$$ where
each $E_n \subset E$ is an ideal and $E = \overline{\cup_n E_n}$.
Note also that each $E_n$ is nuclear since extensions of nuclear
algebras are again nuclear.  Since a locally QD algebra is actually QD
it suffices to show that each $E_n$ is QD.  Since $E_1$ is stably
finite (being a subalgebra of $E$) we have that the boundary map
$\partial : K_1 ( C({\mathbb T}) ) \to K_0 (E_1)$ takes no positive
values.  But then the proof of Proposition~\ref{4.5} shows that if we
assume 4 then $E_1$ will be QD. Proceeding by induction we may assume
that $E_{n - 1}$ is QD.  Since $E_n$ is also stably finite, $E_{n -
1}$ is an ideal in $E_n$ and $E_n/ E_{n - 1} = C({\mathbb T})$,
applying the same argument to the exact sequence $0 \to E_{n - 1} \to
E_n \to C({\mathbb T}) \to 0$ we see that $E_n$ is also QD.

We now show that 3 $\implies$ 5, which will complete the proof.  So
let $A$ be any separable, nuclear, QD $C^*$-algebra and $x \in K_0(A)$
be such that ${\mathbb Z}x \cap K_0^+ (A) = \{ 0 \}$.  Construct a
short exact sequence $0 \to {\cal K}\otimes A \to E \to C({\mathbb T})
\to 0$ such that $\partial (1) = x$.  We will show that $E$ must be
QD.

We can use the $K_0$-Hahn-Banach property to construct an embedding
$\rho : {\cal K}\otimes A \to Q(\oplus_i M_{n_i} ({\mathbb C}))$ such
that $\rho_* (x) = 0$.  Let $D \subset Q(\oplus_i M_{n_i} ({\mathbb
C}))$ be the hereditary subalgebra generated by $\rho({\cal K}\otimes
A)$.  Let $\pi : C({\mathbb T}) \to B(H)$ be any faithful unital
representation such that $\pi (C({\mathbb T})) \cap {\cal K}(H) = \{ 0
\}$.  We first claim that there is an embedding of $E$ into
$(\pi(C({\mathbb T})) + {\cal K}(H)) \otimes \tilde{D}$, where
$\tilde{D}$ is the unitization of $D$.  Indeed, since the embedding
$\rho : {\cal K}\otimes A \to D$ is approximately unital we get a
commutative diagram
$$\begin{CD}
0 @>>> {\cal K}\otimes A @>>> E @>>> C({\mathbb T}) @>>> 0 \\
@. @VVV @VVV @| @. \\
0 @>>> D @>>> F @>>> C({\mathbb T}) @>>> 0,
\end{CD}$$
for some algebra $F$ and the map $E \to F$ is injective.  Since
$\rho_* (x) = 0 \in K_0 (D)$ (by Lemma~\ref{4.4}) and $K_1 (D) = 0$
(by the proof of Lemma~\ref{3.2}) it follows that both boundary maps
arising from the sequence $0 \to D \to F \to C({\mathbb T}) \to 0$ are
zero. Hence we may appeal to the UCT, add on a trivial absorbing
extension and eventually find an embedding of $F$ into $\pi(C({\mathbb
T}))\otimes 1 + {\cal K}(H)\otimes D \subset (\pi(C({\mathbb T})) +
{\cal K}(H)) \otimes \tilde{D}$.

Since $E$ is nuclear it now suffices to show that every nuclear
subalgebra of $(\pi(C({\mathbb T})) + {\cal K}(H)) \otimes \tilde{D}$
is QD.  Hence, by \cite[Prop.~8.3]{Br3} and the Choi-Effros lifting
theorem (\cite{CE}) it suffices to show that there exists a short
exact sequence $$0 \to J \to C \to (\pi(C({\mathbb T})) + {\cal K}(H))
\otimes \tilde{D}\to 0,$$ where $C$ is QD and $J$ contains an
approximate unit consisting of projections which is quasicentral in
$C$ (i.e.\ the extension is quasidiagonal). However, this is now
trivial since $D \subset Q(\oplus_i M_{n_i} ({\mathbb C}))$ implies
that there is a quasidiagonal extension $$0 \to \oplus_i M_{n_i}
({\mathbb C}) \to R \to \tilde{D} \to 0,$$ where $R \subset \Pi_i
M_{n_i} ({\mathbb C})$. But since $X = \pi(C({\mathbb T})) + {\cal
K}(H)$ is nuclear the sequence
$$0 \to (\oplus_i M_{n_i} ({\mathbb C}))\otimes X \to R\otimes X \to
\tilde{D}\otimes X \to 0$$ is exact and since $X$ is unital the
extension is also quasidiagonal.  $\square$

Though Theorem~\ref{4.10} is stated for the class of nuclear QD
$C^*$-algebras a close inspection of the proof shows that this
assumption was only used in the proof of 4 $\implies$ 2. Hence we also
have the following result which applies to individual nuclear
C*-algebras.

\begin{theorem}\label{4.11}
Let $A$ be a separable nuclear QD $C^*$-algebra and consider the following statements.
\begin{enumerate}
\item $A$ has the QD extension property.
\item $A$ has the $K_0$-embedding property.
\item $A$ has the $K_0$-Hahn-Banach property.
\item If $x \in K_0(A)$ is such that ${\mathbb Z}x \cap K_0^+ (A) = \{ 0
\}$ then there exists an embedding $\rho : A \hookrightarrow C$, where
$C$ is QD (but not necessarily separable or nuclear), such that
$\rho_* (x) = 0$.
\item If $x \in K_0(A)$ is such that ${\mathbb Z}x \cap K_0^+ (A) = \{
0 \}$ then there exists a short exact sequence $0 \to {\cal K}\otimes
A \to E \to C({\mathbb T}) \to 0$ where $E$ is QD and $x \in
\partial(K_1 (C({\mathbb T}))) = \partial( {\mathbb Z})$.
\end{enumerate}
Then 1 $\Longleftrightarrow$ 2 $\implies$ 3 $\Longleftrightarrow$ 4
$\Longleftrightarrow$ 5.
\end{theorem}

\textbf{Remark.} There is another version of Theorem~\ref{4.10} where
the class of nuclear C*-algebras is replaced by a class $\mathcal{A}$
of separable C*-algebras with the following closure property. If $0
\to A \otimes \mathcal{K} \to E \to B \to 0$ is exact with $A \in
\mathcal{A}$ and $B$ separable abelian, then $E \in \mathcal{A}$. For
instance $\mathcal{A}$ can be the class of all separable C*-algebras
or the class of all separable exact C*-algebras. Then the statements
1-5 of Theorem~\ref{4.10} formulated for the class $\mathcal{A}$
(rather then for the class of nuclear C*-algebras) are related as
follows: 1 $\Longleftrightarrow$ 2 $\Longleftrightarrow$ 4
$\Longleftrightarrow$ 5 $\implies$ 3.


\begin{thebibliography}{99}


  \bibitem[Bl]{Bl} B. Blackadar, {\it K-theory for operator algebras},
                   Springer-Verlag, New York (1986).

  \bibitem[BK]{BK} B. Blackadar and E. Kirchberg, {\it Generalized inductive
              limits of finite dimensional C*-algebras}, Math. Ann.
              {\bf 307} (1997), 343 - 380.
\bibitem[B]{Brown} L.~G.~ Brown,\emph{ The universal coefficient theorem for ${\rm Ext}$ and quasidiagonality}.
Operator algebras and group representations, Vol. I (Neptun, 1980), 60--64, Monographs
Stud. Math., 17, Pitman, Boston, Mass.-London, 1984.

\bibitem[BD]{BD} L.G. Brown and M. Dadarlat, {\it Extensions of
$C^*$-algebras and quasidiagonality}, J. London Math. Soc. {\bf 53}
(1996), 582 - 600.

\bibitem[BP]{BP} L.G. Brown and G.K. Pedersen, {\it $C^*$-algebras of
real rank zero}, J. Funct. Anal. {\bf 99} (1991), 131 - 149.


\bibitem[Br1]{Br1} N.P. Brown, {\it AF Embeddability of Crossed
  Products of AF algebras by the Integers}, J. Funct. Anal. {\bf 160}
  (1998), 150 - 175.

\bibitem[Br2]{Br2} N.P. Brown, {\it Crossed products of UHF algebras by
some amenable groups}, Hokkaido Math. J. (to appear).

\bibitem[Br3]{Br3} N.P. Brown, {\it On quasidiagonal
$C^*$-algebras}, preprint.


\bibitem[CE]{CE} M.D. Choi and E. Effros, {\it The completely  positive
              lifting problem for C*-algebras}, Ann. of Math. {\bf
              104} (1976), 585 - 609.





  \bibitem[DE1]{DE1} M. Dadarlat and S. Eilers, {\it On the
classification of nuclear $C^*$-algebras}, preprint.

  \bibitem[DE2]{DE2} M. Dadarlat and S. Eilers, {\it Asymptotic
  unitary equivalence in KK-theory}, preprint.

\bibitem[DL]{DL} M. Dadarlat and T.A. Loring, {\it The
K-theory of abelian subalgebras of AF algebras}, J. Reine
Angew. Math. {\bf 432} (1992), 39 - 55.

\bibitem[Da]{Da} K.R. Davidson, {\it C$^{*}$-algebras by Example}, Fields
                   Inst. Monographs vol. 6, Amer. Math. Soc., (1996).

\bibitem[ELP]{ELP} S. Eilers, T.A. Loring and G.K. Pedersen, {\it
Quasidiagonal extensions and AF algebras}, Math. Ann. {\bf 311} (1998),
233 - 249.

\bibitem[Kas]{Kas}
G.G. Kasparov, \emph{The operator ${K}$-functor and extensions of
  {$C^*$}-algebras}, Izv. Akad. Nauk SSSR Ser. Mat. \textbf{44} (1980), no.~3,
  571--636, 719.

\bibitem[EL]{EL} G.A. Elliott and T.A. Loring, {\it AF embeddings of
$C({\mathbb T}^2)$ with a prescribed K-theory}, J. Funct. Anal. {\bf 103} (1992), 1 -
25.

\bibitem[Lo]{Lo} T. Loring, {\it The K-theory of AF embeddings of
the rational rotation algebras}, K-theory {\bf 4} (1991), 227 - 243.



\bibitem[Pa]{Pa} V. Paulsen, {\it Completely bounded maps and dilations},
                 Pitman Research Notes in Mathematics, vol. 146, Longman
                 (1986).


\bibitem[Pe]{Pe1} G. Pedersen, {\it $C^*$-algebras and their
                   automorphism groups}, Academic Press, London (1979).


  \bibitem[Pi]{Pi} M. Pimsner, {\it Embedding some transformation group
                   C$^{*}$-algebras into AF algebras}, Ergod. Th. Dynam.
                   Sys. {\bf 3} (1983), 613 - 626.


\bibitem[PV]{PV} M. Pimsner and D. Voiculescu, {\it Imbedding the
irrational rotation algebras into an AF algebra}, J. Operator Theory
{\bf 4} (1980), 201 - 210.

\bibitem[RS]{RS} J. Rosenberg and C. Schochet, {\it The Kunneth
theorem and the universal coefficient theorem for Kasparov's
generalized $K$-functor}, Duke Math. J. {\bf 55} (1987), 431 - 474.

\bibitem[Sa1]{Sa} N. Salinas, {\it Homotopy invariance of {\rm Ext}(A)}, Duke
Math. J. {\bf 44} (1977), 777 - 794.

\bibitem[Sa2]{Salinas} N.~ Salinas \emph{Relative quasidiagonality and $KK$-theory.} Houston J. Math. \textbf{18}
(1992), no. 1, 97--116.

\bibitem[Sch]{Sch}
C.~Schochet \emph{On the fine structure of the Kasparov groups III: Relative
quasidiagonality}, preprint 1999.

\bibitem[Sp]{Sp} J.S. Spielberg, {\it Embedding C$^{*}$-algebra extensions
                   into AF algebras}, J. Funct. Anal. {\bf 81} (1988),
                   325 - 344.

\bibitem[Vo1]{vo1} D. Voiculescu, {\it A note on quasidiagonal
$C^*$-algebras and homotopy}, Duke Math. J. {\bf 62} (1991), 267 - 271.

\bibitem[Vo2]{vo2} D.~Voiculescu, \emph{Around quasidiagonal operators}
 Integral Equations Operator Theory \textbf{17} (1993),
no. 1, 137--149.



\end{thebibliography}
\end{document}